\documentclass{amsart}
\usepackage{amssymb,latexsym}

\theoremstyle{plain}

\newtheorem{lemma}{Lemma}

\newtheorem{remark}{Remark}

\theoremstyle{definition}
\newtheorem{definition}{Definition}

\theoremstyle{remark}

\numberwithin{equation}{section}

\begin{document}
\title[Lattice Hahn Decomposition ]{Hahn Decomposition Theorem \\ of Signed Lattice Measure}

\author{Jun Tanaka}
\address{University of California, Riverside, USA}
\email{juntanaka@math.ucr.edu, yonigeninnin@gmail.com}

\date{January, 15, 2009}

\begin{abstract}
  In this paper, we will define a signed Lattice measure on $\sigma$-algebras, as well as give the definition of positive and negative Lattice. Herein, we will show that the Hahn Decomposition Theorem decomposes any space X into a positive lattice A and a negative Lattice B such that $A \vee B$ =X and the signed Lattice measure of $A \wedge B $ is 0.
\end{abstract}

\maketitle

\section{\textbf{Introduction}}
In this paper, we will show Hahn Decomposition Theorem of Signed Lattice Measure when the Lattice do not necessarily satisfy the law of non-contradiction. We consider a set X and a Lattice $\mathcal{L}$ of subsets of X, which may not satisfy the law of non-contradiction and which contain $\emptyset$ and X.

We will define a signed Lattice measure on $\sigma$-algebras, and we will show that the Lattice Hahn Decomposition Theorem decomposes any set X into a positive Lattice A and a negative Lattice B such that $A \vee B $=X and the signed Lattice measure of $A \wedge B $ is 0.

This theory discussed in this paper is a Lattice version of the paper \cite{Junhahndeco}. In the paper \cite{Junhahndeco}, J. Tanaka showed Hahn Decomposition Theorem of Signed Fuzzy Measure. As one already saw, Fuzzy measure is a classical measure, provided that fuzzy sets are restricted to classical sets. As in classical measure theory, he defined a fuzzy signed measure on $\sigma$-algebras, and he showed that the Fuzzy Hahn Decomposition Theorem decomposes any space X into a positive set A and a negative set B such that A+B=X and the signed measure of $A \wedge B $ is 0.

 Unfortunately, the monotonicity prohibits the signed Lattice measure from being a signed measure in a classical sense. Please note that monotonicity is not required in the definition of classical signed measure. In Section \ref{sec:2}, we will define a signed Lattice measure on $\sigma$-algebras, as well as give the definition of positive and negative Lattice. Furthermore, we will show that any countable union of positive Lattice is a positive Lattice, as well as the following proposition: if E is a Lattice such that 0 $< \nu (E)   <  \infty$, then there is a positive Lattice A $\leq E$ with $ \nu (A) >$ 0. Finally, we will show the Hahn Decomposition Theorem for signed Lattice measure in the section titled Main Result.

\section{\textbf{Preliminaries : Extension of Lattice}}\label{modernS:1}

In this section, we shall briefly review the well-known facts about lattice theory (e.g. Birkhoff \cite{Birk} ), propose an extension lattice, and investigate its properties. (L,$\wedge  $,$\vee  $) is called a lattice if it is closed under operations $\wedge  $ and $\vee  $ and satisfies, for any elements x,y,z in L:

(L1) the commutative law: x $\wedge$ y = y $\wedge$ x and x $\vee$ y = y $\vee$ x

(L2) the associative law:
\[
\begin{aligned}
x \wedge ( y \wedge z ) = (x \wedge y) \wedge z \ \ \text{and} \ \   x \vee ( y \vee z ) = (x \vee y) \vee z
\end{aligned}
\]

(L3) the absorption law:    x $\vee$ ( y $\wedge$ x ) =x  and   x $\wedge$ ( y $\vee$ x ) = x.

Hereinafter, the lattice (L,$\wedge  $,$\vee  $) will often be written as L for simplicity.

A mapping h from a lattice L to another $L'$ is called a lattice-homomorphism, if it satisfies
\[
\begin{aligned}
h ( x \wedge y ) = h(x) \wedge h(y)  \ \ \text{and} \ \   h ( x \vee y ) = h(x) \vee h(y) , \forall x,y \in L.
\end{aligned}
\]

If h is a bijection, that is, if h is one-to-one and onto, it is called a lattice-isomorphism; and in this case, $L'$ is said to be lattice-isomorphic to L.

A lattice (L,$\wedge  $,$\vee  $) is called distributive if, for any x,y,z in L,

(L4) the distributive law holds:
\[
\begin{aligned}
x \vee ( y \wedge z ) = (x \vee y) \wedge (y \vee z)  \ \ \text{and}  \  \   x \wedge ( y \vee z ) = (x \wedge y) \vee ( y \wedge z)
\end{aligned}
\]

A lattice L is called complete if, for any subset A of L, L contains the supremum $\vee$ A and the infimum $\wedge$ A. If L is complete, then L itself includes the maximum and minimum elements, which are often denoted by 1 and 0, or I and O, respectively.

A distributive lattice is called a Boolean lattice if for any element x in L, there exists a unique complement $x^C$ such that
\[
\begin{aligned}
& x \vee x^C = 1  \ \ \ \ \ \ \ \ \text{(L5) the law of excluded middle}  \\
\ \ &x \wedge x^C = 0 \ \ \ \ \ \ \ \ \text{(L6) the law of non-contradiction}.
\end{aligned}
\]

Let L be a lattice and $\cdot^{c}$: L $\rightarrow $ L be an operator. Then $\cdot^{c}$ is called a lattice complement in L if the following conditions are satisfied:
\[
\begin{aligned}
&\text{(L5) and (L6)};                     & \forall x \in L, \ x^C \vee x = I \ \text{and} \ x^C \wedge x = 0,   \\
&\text{(L7) the law of contrapositive};    &  \forall x, y \in L, x \leq y \Rightarrow x^C \geq y^C ,  \\
&\text{(L8) the law of double negation};   &  \forall x \in L, (x^C)^C = x    .
\end{aligned}
\]
\begin{definition} {\textbf{Complete Heyting algebra (cHa)}}

A complete lattice is called a complete Heyting algebra (cHa) if
\[
\begin{aligned}
\vee_{i \in I} \ ( x_{i} \wedge y ) =  (\vee_{i \in I} \  x_{i} ) \wedge y
\end{aligned}
\]
holds for $\forall x_{i} , y \in L $ ($i \in I$); where I is an index set of arbitrary cardinal number.
\end{definition}

It is well-known that for a set E, the power set P(E) = $2^E$. The set of all subsets of E is a Boolean algebra.

\section{\textbf{Definitions and Lemmas}}\label{sec:2}

Throughout this paper, we will consider Lattices as complete Lattices which obey (L1)-(L8) except for (L6) the law of non-contradiction.

\begin{definition}

Unless otherwise stated, X is the entire set and $\mathcal{L} $ is a Lattice of any subset sets of X. If a Lattice $\mathcal{L} $ satisfies the following conditions, then it is called a Lattice $\sigma$-algebra;

(1) $\forall  h \in \mathcal{L}$, $ h^C \in \mathcal{L}$

(2)  if $h_{n} \in \mathcal{L}$ for n=1,2,3...., then $ \vee^{\infty} h_{n} \in \mathcal{L} $.

We denote $\sigma ( \mathcal{L} )$ as the Lattice $\sigma$-Algebra generated by $\mathcal{L}$.

\end{definition}

\begin{definition}\label{De:2}
If m : $\sigma ( \mathcal{L} )$ $\mapsto$ $\mathbb{R} \cup \{\infty \}$ satisfies the following properties, then m is called a Lattice measure on the Lattice $\sigma$-Algebra $\sigma ( \mathcal{L} )$.

(1) m $ (\emptyset) $ = m $ (0)$ = 0.

(2) $\forall h , g  \in \sigma ( \mathcal{L} )$ s.t. $m( h) , m(g) \geq   $ 0 : $h \leq g$ $ \Rightarrow  $  m $ ( h) \leq $ m $(g)$.

(3) $\forall h , g  \in \sigma ( \mathcal{L} )$ : m$( h \vee g ) + $ m$( h \wedge g ) = $ m$( h)  +  $ m$ (g)   $.

(4) if $ h_{n}  $ $\subset  \sigma ( \mathcal{L} ), n \in N$ such that $h_{1} \leq   h_{2} \leq    \cdots  \leq  h_{n} \leq  \cdots $,  then $m ( \vee^{\infty} h_{n} ) =  \lim $ m$( h_{n} )$.

\end{definition}

Let $m_{1}$ and $m_{2}$ be Lattice measures defined on the same Lattice $\sigma$-algebra $\sigma ( \mathcal{L} )$. If one of them is finite, the set function $m ( E ) =  m_{1}( E )  - m_{2}( E )$ , $E \in \sigma ( \mathcal{L} )$ is well defined and countably additive on $\sigma ( \mathcal{L} )$. However, it is not necessarily nonnegative; it is called a signed Lattice measure.

\begin{definition}\label{De:1}
 By a signed Lattice measure on the measurable Lattice (X, $\sigma ( \mathcal{L} )$) we mean $\nu$ : $\sigma ( \mathcal{L} )$ $\mapsto$ $\mathbb{R} \cup \{\infty \}$ or $\mathbb{R} \cup \{- \infty \}$, satisfying the following property:

(1) $\nu (\emptyset) $ = $\nu (0)$ = 0.

(2) $\forall h , g  \in \sigma ( \mathcal{L} )$ s.t. $\nu( h) , \nu(g) \geq   $ 0 : $h \leq g$ $ \Rightarrow  $  $ \nu( h) \leq \nu(g)$.

\ \ \ \ \ $\forall h , g  \in \sigma ( \mathcal{L} )$ s.t. $\nu( h) , \nu(g) \leq   $ 0 : $h \leq g$ $ \Rightarrow  $  $  \nu(g)  \leq  \nu( h)$.

(3) $\forall h , g  \in \sigma ( \mathcal{L} )$ : $\nu( h \vee g ) +   \nu( h \wedge g ) = \nu( h)  +   \nu(g)   $.

(4) if $ h_{n}  $ $\subset  \sigma ( \mathcal{L} ), n \in N$ such that $h_{1} \leq   h_{2} \leq    \cdots  \leq  h_{n} \leq  \cdots $,  then $m ( \vee^{\infty} h_{n} ) =  \lim $ m$( h_{n} )$.

This is meant in the sense that if the left-hand side is finite, the limit on the right-hand side is convergent, and if the left-hand side is $ \pm \infty $, then the limit on the right-hand side diverges accordingly.
\end{definition}

\begin{remark}
The signed Lattice measure is a Lattice measure when it takes only positive values. Thus, the signed Lattice measure is a generalization of Lattice measure.

\end{remark}

\begin{definition}
A is a positive Lattice if for any Lattice measurable set E in A, $\nu (E) \geq 0$. Similarly, B is a negative Lattice if for any Lattice measurable set E in B, $\nu (E) \leq 0$.

\end{definition}

\begin{lemma}\label{Le:1}
Every sublattice of a positive Lattice is a positive Lattice and any countable union of positive Lattices is a positive Lattice.

\begin{proof}
The first claim is clear. Before we show the second claim, we need to show that every union of positive Lattice is a positive Lattice. Let A, B be positive Lattices and E $\leq$ A $\vee $ B be a measurable Lattice. By (2) in Definition $\ref{De:1}$, 0 $\leq  $ $\nu (B \wedge E) $ - $\nu (A \wedge B \wedge E) $. By (3), $\nu (E) \geq$ 0. Now by induction, every finite union of positive Lattice is a positive Lattice. Let $A_{n}$ be a positive Lattice for all n and E  $\leq$ $\vee $ $A_{n}$ be a measurable Lattice. Then $E_{m}$ := $ E \wedge \vee_{n=1}^{m} A_{n}  $ =  $\vee_{n=1}^{m}  E  \wedge A_{n}  $. Then $E_{m}$ is a measurable Lattice and a positive Lattice. In particular,  $E_{m} \leq E_{m+1}$ for all n and $ E$ = $ \vee_{n=1}^{\infty} E_{m} $. Thus 0 $\leq  \lim \nu ( E_{m} ) $ = $\nu ( E ) $. Therefore $\vee $ $A_{n}$ is a positive Lattice.

\end{proof}
\end{lemma}

\begin{lemma}\label{Le:2}

Let E be a  measurable Lattice such that 0 $< \nu (E)   <  \infty$. Then there is a positive Lattice A $\leq E$ with $ \nu (A) >$ 0.

\begin{proof}
If E is a positive Lattice, we take A=E. Otherwise, E contains a Lattice of negative measure. Let $n_{1}$ be the smallest positive integer such that there is a measurable Lattice $E_{1}  \subset $ E with $ \nu ( E_{1} ) < - \frac{1}{n_{1}} $. Proceeding inductively, if E$\wedge \wedge_{j=1}^{k-1} E_{j}^C  $ is not already a positive Lattice, let $n_{k}$ be the smallest positive integer for which there is a measurable Lattice $ E_{k} $ such that $ E_{k} \leq E \wedge \wedge_{j=1}^{k-1} E_{j} $ and $\nu ( E_{k}) < - \frac{1}{n_{k}} $.

Let A = $( \vee  E_{k})^C $.

Then $\nu (E)$ = $ \nu ( E \wedge A  ) + \nu(E   \wedge   \vee  E_{k}  )  $ = $ \nu (E \wedge  A  ) + \nu(   \vee  E_{k}  )  $. Since  $\nu (E)$ is finite, $ \lim_{n \rightarrow \infty}   \nu(   \vee^{n} E_{k}  )$ is finite and $ \nu(   \vee  E_{k}  ) \leq$ 0. Since $\nu (E) >$ 0 and $ \nu(   \vee  E_{k}  ) \leq$ 0, $ \nu (E \wedge A  ) >$ 0.

We will show that A is a positive Lattice. Let $  \epsilon  >$ 0. Since $ \frac{1}{n_{k}} \rightarrow  $ 0, we may choose k such that $-\frac{1}{n_{k} - 1} \  > \  - \epsilon $. Thus A contains no measurable Lattice of measure less than $ - \epsilon $. Since $ \epsilon $ was an arbitrary positive number, it follows that A can contain no Lattice of negative measure and so must be a positive Lattice.

\end{proof}
\end{lemma}

\section{\textbf{Main result: Lattice Hahn Decomposition}}\label{sec:4}

Without loss of generality, let's omit + $\infty $ as a value of $\nu$. Let $\lambda$ = $\sup \{ \nu (A) : A $ is a positive Lattice $      \}   $.

Then $\lambda   \geq$ 0 since $\nu (\emptyset) $ = 0.

Let {$A_{i} $} be a sequence of positive Lattices such that $\lambda$ = $\lim  \nu (A_{i} )$ and $  A$ = $\vee A_{i} $. By Theorem $\ref{Le:1}$, A is a positive Lattice and $\lambda \geq \nu (A)$.

$\vee^{n}  A_{i} \leq $ $A$ for any n implies  $\nu ( \vee^{n}  A_{i}) \geq$ 0 for any n. Thus $\lambda$ = $\lim \nu (   A_{i} )$ = $ \nu (A) $ = 0.

Let $E \leq A^C$ be a positive Lattice. Then $\nu (E) \geq$ 0 and $A \vee E $ is a positive Lattice. Thus $\lambda   \geq$ $ \nu ( A \vee E    ) = \nu ( A ) + \nu ( E )    - \nu ( A  \wedge E  )   $ = $\lambda + \nu ( E) - \nu ( A  \wedge E  )    $. Thus $\nu (E) = \nu ( A  \wedge E  )    $. We have that $\nu (E)$ = 0 since $ A  \wedge E \leq A \wedge A^C$ and $\nu ( A \wedge A^C ) $ = 0.

Thus, $A^C$ contains no positive sublattice of positive measure and hence no sublattice of positive measure by Lemma $\ref{Le:2}$. Consequently, $A^C$ is a negative Lattice.

\section{\textbf{Conclusion}}
Let X be an entire set. Then by the previous theorem, we find such a positive Lattice A and a negative Lattice B (= $A^C$). By the Lattice measurability of $\nu$, $\nu ( A  \wedge  A^C    ) $ = 0. $A \vee A^C$ = X. These characteristics provide the following: X = $ A \cup  B$ and $ A \cap B = \emptyset$ in the classical set sense.



\end{document}